\documentclass[reqno,10pt]{amsart}
\usepackage[english]{babel}
\usepackage{amsthm}
\usepackage{amsmath}
\usepackage{amssymb}
\usepackage{amsfonts}
\usepackage{mathrsfs}
\usepackage{todonotes}
\usepackage{comment}

\usepackage[utf8]{inputenc}
\usepackage{latexsym}
\usepackage{verbatim}

\usepackage[all]{xy}
\usepackage{dsfont}
\usepackage[bookmarksnumbered,colorlinks,urlcolor=blue]{hyperref}
\usepackage{graphicx}
\usepackage{graphics}
\usepackage{float}
\usepackage{enumerate}
\def\bb#1\eb{\textcolor{blue}
{#1}} %
\def\br#1\er{\textcolor{red}
{#1}} %
\newcommand{\R}{\mathds R}

   \def\br#1\er{\textcolor{red}{#1}} %
      \def\bb#1\eb{\textcolor{blue}{#1}} %

\title[Precompact Holonomy]{Compact affine manifolds with precompact holonomy are geodesically complete}

\author[L. Ak\'e]{Luis Alberto Ak\'e Hau}
\address{Departamento de \'Algebra, Geometr\'{\i}a y Topolog\'{\i}a,  Universidad de M\'alaga
\hfill\break\indent
Facultad de Ciencias, Campus Universitario de Teatinos,
\hfill\break\indent 29080 M\'alaga, Spain}
\email{luisake@uma.es}

\author[M. S\'anchez]{Miguel S\'anchez}
\address{Departamento de Geometr\'{\i}a y Topolog\'{\i}a,  Universidad de Granada 
\hfill\break\indent
Facultad de Ciencias, Campus Fuentenueva s/n,
 \hfill\break\indent E-18071 Granada, Spain}
\email{sanchezm@ugr.es}

\begin{document}
\newtheorem{thm}{Theorem}[section]
\newtheorem{prop}[thm]{Proposition}
\newtheorem{lemma}[thm]{Lemma}
\newtheorem{cor}[thm]{Corollary}
\theoremstyle{definition}
\newtheorem{defi}[thm]{Definition}
\newtheorem{notation}[thm]{Notation}
\newtheorem{exe}[thm]{Example}
\newtheorem{conj}[thm]{Conjecture}
\newtheorem{prob}[thm]{Problem}
\newtheorem{rem}[thm]{Remark}
\newtheorem{conv}[thm]{Convention}
\newtheorem{crit}[thm]{Criterion}
\newtheorem{claim}[thm]{Claim}

\newcommand{\ben}{\begin{enumerate}}
\newcommand{\een}{\end{enumerate}}

\newcommand{\bit}{\begin{itemize}}
\newcommand{\eit}{\end{itemize}}

\begin{abstract}
This note proves  the geodesic completeness of any compact  manifold endowed with a linear connection such that the closure of its holonomy group is compact. 
\end{abstract}

\maketitle 

\vspace*{-.5cm}

\section{Introduction}\label{section1}

The purpose of the present note is to prove the following result:

\begin{thm}\label{t}
Let $M$ be a (Hausdorff, connected,  smooth) compact $m$-manifold endowed with a linear connection $\nabla$ and let $p\in M$. If the holonomy group  Hol$_{p}(\bigtriangledown)$ (regarded as a subgroup of the group Gl($T_pM$) of all the linear automorphisms of the tangent space at $p$, $T_pM$) has  compact closure, then  $(M,\bigtriangledown)$ is geodesically complete.
\end{thm}
Some comments on the completeness of compact affine manifolds are in order. There are several results when $\nabla$ is {\em flat} and, therefore, there exists
an atlas $\mathcal{A}$ whose transition functions are affine maps of $\R^m$; in this case, the linear parts will lie in some subgroup $G$ of the (real) general linear group Gl$(m)$. In fact, a well-known conjecture by Markus states that compact affine flat manifolds which are  unimodular (i.e.,  $G$ can be chosen in the special linear group Sl$(m)$) must be complete. Carri\`ere \cite{Ca} introduced an invariant for linear groups, the {\em discompacity}, which measures the non-compactness of the group by analyzing the degeneration of images of the unit sphere under the action of sequences of its elements. When the closure $\bar G$ is compact, the discompacity is equal to 0; Markus conjecture was proven in  \cite{Ca} under the assumption that the discompactness of $G$ is at most 1.
Other results on the structure of unimodular manifolds (see for example, \cite{Ca2}) can be also regarded as partial answers to that conjecture.

When $\nabla$ is the Levi-Civita connection of a Riemannian metric $g$, the (geodesic) completeness of $\nabla$, which follows directly from Hopf-Rinow theorem, can be reobtained from Theorem \ref{t}; indeed, Hol$_{p}(\bigtriangledown)$ becomes a subgroup of the (orthogonal) group of linear isometries O$(T_pM)$, which is compact. However, when $g$ is an indefinite semi-Riemannian metric of index $\nu$ ($0<\nu<m$),
the corresponding orthogonal group O$_\nu (T_pM)$ is non-compact and completeness may not hold; 
the Clifton-Pohl torus  (see for instance \cite[Example 7.16]{ON}) is   
a well-known example.
There are some results  that assure the completeness of a compact semi-Riemannian manifold, among them either to admit $\nu$ pointwise independent conformal Killing vector fields  which span a negative definite subbundle of $TM$   \cite{RS0, MSanchez2}, or to be homogeneous \cite{Marsden} (local homogeneity is also enough in dimension 3 \cite{Z}, and to be conformal to a homogeneous manifold is enough in any dimension  \cite{MSanchez2}; see \cite{MSanchez} for a review). 

In the particular case that the semi-Riemannian manifold $(M,g)$ is Lorentzian ($\nu=1<m$), compactness implies completeness in other relevant cases, such as when $g$ is flat. Indeed, taking if necessary a  finite covering, this is a particular case   of Markus' conjecture where $G$ can be regarded as the restricted Lorentz group SO$^{\uparrow}_1(m)$ (i.e., the connected component of the identity of the Lorentz group O$_1(m)$) and, as also proven by 
Carri\`ere in \cite{Ca}, the discompactness  of SO$^{\uparrow}_1(m)$ is equal to 1. It is worth pointing out that, when the group $G$ determined by a compact flat affine manifold lies in the  group of Lorentzian similarities (generated by homotheties and O$_1(m)$) but not in the Lorentz group, then the connection is incomplete \cite{A}. Moreover, Klinger \cite{Kl} extended Carri\`ere's result by showing that any compact Lorentzian manifold of constant curvature is complete, and  Leistner and Schliebner 
 \cite{AbeH} proved that completeness also holds in the case of Abelian holonomy (compact pp-waves). 

These  semi-Riemannian  results are independent of Theorem \ref{t}; in fact, Guti\'errez and M\"uller  \cite{Guti} have proven recently that, for a Lorentzian metric $g$, the compactness of $\overline{\hbox{Hol}(g)}$ 
implies the existence of a timelike parallel vector field in a finite covering. 
This conclusion (combined with the cited one in \cite{RS0}) also gives an alternative proof of Theorem \ref{t} in the particular case that $\nabla$ comes from a Lorentzian metric. 
In any case, the proof of our theorem is very simple and extends or complements the previous results.

\section{Proof of  Theorem \ref{t}} \label{s2}

Assume that there exists 
an incomplete geodesic  $\gamma:[0,b) \rightarrow M$, $b < \infty$. By using the  compactness of $M$, choose any sequence $\{t_{n}\}_n\nearrow b$  such that $\{\gamma(t_{n})\}_n$ converges to some $p \in M$. It is well-known then that the sequence of velocities $\{\gamma'(t_{n})\}_n$ cannot converge in $TM$ as $\gamma'$ is the integral curve of the geodesic vector field on $TM$ (see for example Prop. 3.28 and Lemma 1.56 in \cite{ON} or \cite[Section~3]{MSanchez}).
Consider a  normal (starshaped) neighbourhood $U$ of $p$ (see for example \cite{p, white} for background results on linear connections). With no loss of generality, we will assume that  
$\{\gamma(t_{n})\}_n\subset U$ and will arrive at a contradiction with the compactness of $\overline{\hbox{Hol}_p(\nabla)}$.

Consider  the loops at $p$ given by 
$\alpha_{n}=\rho_{n}^{-1} \star \gamma_{[t_{1},t_{n}]} \star \rho_{1}$, where $\star$ denotes the concatenation of the corresponding curves and $\rho_{n}: [0,1]\rightarrow U$ is the radial geodesic from $p$ to $\gamma(t_{n})$
for all $n=1,2 \dots$ Put $v_p=\tau_{\rho_{1}^{-1}}(\gamma'(t_{1})) \in T_{p}M$. For any curve $\alpha: [a,b]\rightarrow$, let $\tau_\alpha$ be  
the parallel transport between its endpoints. As $\gamma$ is a geodesic: 
$$v_n:= \tau_{\alpha_{n}}(v_{p})=\tau_{\rho_{n}^{-1}} \circ \tau_{\gamma_{[t_{1},t_{n}]}} \circ 
\tau_{\rho_{1}}(\tau_{\rho_{1}^{-1}}(\gamma'(t_{1})))=\tau_{\rho_{n}^{-1}}(\gamma'(t_{n}))$$
(in particular, $v_1=v_p$). The compactness of $\overline{\hbox{Hol}_p(\nabla)}$ implies that the sequence $\{v_n\}_n$ is contained in a compact subset  $K\subset T_pM$. So, it is enough to check that this property also implies that  $\{\gamma'(t_{n})\}_n$ is included in a compact subset of $TM$.

With this aim, let $\tilde K\subset \exp^{-1}(U)$ be a starshaped compact neighborhood of $0\in T_pM$ and, for each $u\in \tilde K$, let $\rho_u: [0,1]\rightarrow U$ be the radial geodesic segment with initial velocity $u$. As the map $$\xi: \tilde K\times K\rightarrow TM, \qquad (u,v)\mapsto \xi(u,v)=\tau_{\rho_u}(v)$$ 
is continuous, its image $\xi(\tilde K,K)$ is compact. This set contains all $\{\gamma'(t_n)\}_n$ up to a finite number (recall that $\gamma(t_n)$ lies in $\exp_p(\tilde K)$ and $\rho_{u_n}=\rho_n$ for $u_n=\exp^{-1}(\gamma(t_n))$, i.e.,  $\gamma'(t_n)=\xi(u_n,v_n)$), as required.   

\section*{Acknowledgments}
The support of the  projects
MTM2013-47828-C2-1-P and  MTM2013-47828-C2-2-P (Spanish MINECO with FEDER funds) is acknowledged. The first-named author also acknowledges a 
grant funded by the Consejo Nacional de Ciencia y Tecnolog\'ia (CONACyT), M\'exico.

\end{document}